\documentclass[12pt,oneside]{amsart}
\usepackage{tkz-euclide}
\usepackage{tkz-euclide}
\usetikzlibrary{positioning}
\usepackage{geometry}                % See geometry.pdf to learn the layout options. There are lots.
\usepackage{graphicx}
\usepackage{amssymb}

\usepackage{amsmath,amsthm,amscd,amssymb}
\usepackage{latexsym}
\usepackage[colorlinks,citecolor=red,pagebackref,hypertexnames=false]{hyperref}
\usepackage{ulem}
\usepackage{soul}
\numberwithin{equation}{section}
\theoremstyle{plain}
\newtheorem{theorem}{Theorem}[section]
\newtheorem{lemma}[theorem]{Lemma}
\newtheorem{corollary}[theorem]{Corollary}

\theoremstyle{definition}
\newtheorem{definition}[theorem]{Definition}

\newtheorem{example}[theorem]{Example}

\theoremstyle{remark}
\newtheorem{remark}[theorem]{Remark}

\newtheorem{case[theorem]}{Case}

\def\be{\begin{equation}}
\def\ee{\end{equation}}

\def\new{}

\def\R{\mathbb R}
\def\hd{{\dim_{\mathcal H}}}

\date{March 14, 2025}      
\author{A. Greenleaf, A. Iosevich, and K. Taylor}
\address{Department of Mathematics, University of Rochester, Rochester, NY}
\email{allan@math.rochester.edu}
\address{Department of Mathematics, University of Rochester, Rochester, NY}
\email{iosevich@gmail.com}
\address{Department of Mathematics, The Ohio State, Columbus, OH}
\email{taylor.2952@osu.edu}
\thanks{The first listed author was supported in part by the National Science Foundation under grant NSF DMS-2204943. The second listed author was supported in part by the National Science Foundation under grant NSF DMS - 2154232. The third listed author is supported in part by the Simons Foundation Grant no. 523555.}

\begin{document}

\title{Realizing trees of configurations in thin sets}

\maketitle

\begin{abstract} Let $\phi(x,y)$ be a continuous function, smooth away from the diagonal, such that,  for some $\alpha>0$, 
the associated generalized Radon transforms
\begin{equation} \label{Radon} R_t^{\phi}f(x)=\int_{\phi(x,y)=t} f(y) \psi(y) d\sigma_{x,t}(y) \end{equation} 
map $L^2({\mathbb R}^d) \to L^2_{\alpha}({\mathbb R}^d)$ for all $t>0$. 
Let $E$ be a compact subset of ${\mathbb R}^d$ for some $d \ge 2$, and suppose that the 
Hausdorff dimension of $E$ is $>d-\alpha$. 
We show that any tree graph $T$ on $k+1$ ($k \ge 1$) vertices is \new{stably} realizable in $E$, 
in the sense that \new{for each $t$ in some open interval} there exist distinct 
$x^1, x^2, \dots, x^{k+1} \in E$ 
%and $t>0$ 
such that the $\phi$-distance $\phi(x^i, x^j)=t$ 
for all pairs $(i,j)$  corresponding to the edges of  $T$. 

\new{We  extend this result to  trees whose edges are prescribed by more complicated point configurations, 
such as congruence classes of triangles.}

\iffalse 
 where points $x^1, x^2, \dots, x^{k+1}$ are 
connected by a hyper-edge if $\phi(x^i,x^j)=t_{ij}$ 
for the set of pairs $(i,j)$ in the edge set of a graph $H$, 
and $\{t_{ij}\}$ is a pre-assigned collection of positive distances. 
A chain of hyperedges is a sequence 
$$\{(x^1_1, \dots, x^{k+1}_1), (x^1_2, \dots, x^{k+1}_2), \dots, (x^1_n,\dots, x^{k+1}_n)\},$$ 
where $(x^1_i, \dots, x^{k+1}_i)$ and $(x^1_{i+1}, \dots, x^{k+1}_{i+1})$, $1 \leq i \leq n$, share one and only one vertex $x^J_i$, with 
$J$ pre-assigned. \edit
A hyper-graph tree is defined analogously. The key estimate behind these arguments is a variant of (\ref{Radon}) 
where the average over the level set of $\phi(x,\cdot)$ is replaced by the average over an intersection of such surfaces. In the process, 
we develop a general graph theoretic paradigm that reduces a variety of point configuration questions to operator theoretic estimates. 

\fi

\end{abstract}  

\maketitle

\section{Introduction} 

The celebrated Falconer distance conjecture (see, e.g., \cite{Fal86,M95,Mat15}) states that if the 
Hausdorff dimension of a compact set $E \subset {\Bbb R}^d$, $d \ge 2$, is greater than $\frac{d}{2}$, 
then the Lebesgue measure of the distance set $\Delta(E)=\{|x-y|: x,y \in E \}$ is positive. 
Until recently, the best results  known were due to Wolff  \cite{W04} in two dimensions and 
Erdo\~{g}an \cite{Erd05} in higher dimensions. They proved that Lebesgue measure of the distance set 
is positive if the Hausdorff dimension of $E$ satisfies $\hd(E)>\frac{d}{2}+\frac{1}{3}$. When $d=2$, 
Orponen (\cite{O15}) proved that, under the additional assumption that $E \subset {\Bbb R}^2$ is 
Ahlfors-David regular, if $\hd(E)>1$, then the packing dimension 
of  $\Delta(E)$ is $1$.

Currently,  the best known exponent threshold for the Falconer distance problem in two dimensions is 
$\frac{5}{4}$, due to Guth, the second listed author, Ou and Wang (\cite{GIOW20}). In higher 
dimensions, the best exponent in odd dimensions, recently established by 
Du, Ou, Ren and Zhang (\cite{DORZ23}), is $\frac{d}{2}+\frac{1}{4}-\frac{1}{8d+4}$;
see \cite{DIOWZ21} for even dimensions. 

The Falconer problem has many variations, where distance is replaced by more general $k$-point configurations,
which need not be scalar-valued.
For $p\in\mathbb N$ and $k\ge 2$, let $\Phi:\left(\R^d\right)^k \to \R^p$ be a 
continuous function which is smooth (except possibly on a lower dimensional set).
For a compact $E\subset\R^d$, the {\it $\Phi$-configuration set} of $E$ \cite{GGIP12} is the compact set
$$\Delta_\Phi(E):=\left\{\, \Phi\left(x^1,\dots,x^k\right): \, x^1,\dots,x^k\in E\, \right\}\subset \R^p,$$
and one can look for lower bounds on $\hd(E)$ ensuring that $\Delta_\Phi(E)$ has positive Lebesgue measure in $\R^p$.

A further  variation on the Falconer problem,
originating in the result of Mattila and Sj\"olin \cite{MS99} for the distance set,
 seeks to determine values $s_\Phi$ so that   $\hd(E)>s_\Phi$  guarantees that 
$\Delta_\Phi(E)$ has {\it nonempty interior} in $\R^p$,
in which case $\Phi$ is said to be a {\it Mattila-Sj\"olin function}. 
See \cite{MS99,IMT11,GIT19,GIT20,PRA21,KPS21,GIT22,PRA22,GGPP23,ITU16} for results of this type. 
\medskip

A particularly interesting example arises when the Euclidean distance $|x-y|$ is replaced 
by a more general function $\phi(x,y)$. 
For a compact  $E\subseteq \mathbb{R}^d$, for some $d\geq 2$,  
and a $\phi: {\Bbb R}^d \times {\Bbb R}^d \to {\Bbb R}$, continuous and smooth away from the diagonal,
we define the {\it generalized distance set} 
\begin{equation} \label{gendistdef} \Delta_{\phi}(E)=\{\phi(x,y): x,y \in E\}\subset \R. \end{equation}
Eswarathasan and the second and third listed authors proved in \cite{EIT11} that if $\phi$ satisfies the non-vanishing Monge-Amp\`ere 
determinant condition, 
\begin{equation} det
\begin{pmatrix}
 0 & \nabla_{x}\phi \\
 -{(\nabla_{y}\phi)}^{T} & \frac{\partial^2 \phi}{dx_i dy_j}
\end{pmatrix}
\neq 0,
\end{equation} on the set $\{(x,y): \phi(x,y)=t \}$, and if $E \subset {\Bbb R}^d$, $d \ge 2$, 
is a compact set with $\hd(E)>\frac{d+1}{2}$, then the Lebesgue measure of $\Delta_{\phi}(E)$ is positive.
A particularly compelling case arises when $E$ is a subset of a compact Riemannian manifold 
without  boundary or conjugate points and $\phi$ is the induced distance function.  
The second listed author, Liu and Xi proved in two dimensions (\cite{ILX2020}) that if $\hd(E)>5/4$
%$E$, a subset of a compact two-dimensional Riemannian manifold $M$ without  
%boundary, has Hausdorff dimension $>\frac{5}{4}$, and $\phi$ is a Riemannian metric on $M$,
 then the Lebesgue measure of $\Delta_{\phi}(E)$ is positive, 
matching the exponent obtained in the Euclidean case  of \cite{GIOW20}.

\vskip.125in 

The main thrust of this paper is to develop a  general technique to study finite point configurations of a graph-theoretic 
nature in Euclidean space and Riemannian manifolds, 
and apply it to resolve several open problems. 
We let  $G=(\mathbb V,\mathbb E)$ denote an undirected graph on 
$k$ vertices. The {\it edge map} of $G$ is ${\mathcal E}_G:\mathbb V\times \mathbb V\to \{0,1\}$, ${\mathcal E}_G(i,j)=1$ if $i \not=j$ 
and  the $i$'th and $j$'th vertices are connected by an 
edge in $\mathbb E$  , and $0$ otherwise. 
\bigskip

\begin{definition} ({\bf Generalized Distance Graph}) \label{graphdefinition}
A continuous $\phi: {\Bbb R}^d \times {\Bbb R}^d \to {\Bbb R}$, smooth away from the diagonal and such that 
$\phi(x,y)=\phi(y,x)$,  a compact
$E\subset{\Bbb R}^d$, and a $t>0$,
define the {\it generalized distance graph} $G_{\phi,t}(E)$, 
whose vertices  are the points in $E$, and for which two vertices $x,y\in E$, 
$x \not=y$, are connected by an edge iff $\phi(x,y)=t$.

We  say that an (abstract)  connected finite graph $G$ can be {\it realized} in $E$ if there exists $t>0$ such that $G$ 
is isomorphic to a subgraph of $G_{\phi,t}(E)$; furthermore, $G$ is said to be {\it stably realized} in $E$ if the set of such $t$ 
has nonempty interior. 
 \end{definition}

%\vskip.125in 

Bennett and the second and third listed authors of this paper proved (\cite{BIT15}) that if 
$E \subset {\Bbb R}^d$, $d \ge 2$, 
with $\hd(E)>\frac{d+1}{2}$, $\phi(x,y)=|x-y|$, and $G$ is a {\it path} (or {\it chain}), then $G$ can be stably realized in $E$. 
In \cite{IT19}, the second and  third listed authors extended this result %, in the case when $\phi(x,y)=|x-y|$, 
to the more general case when $G$ is a {\it tree}. % graphs. 
\smallskip

We note that trees and chains have also been considered for notions of size other than Hausdorff dimension, see for instance \cite{McTconstantGap, McTchains}, where McDonald and the third listed author prove that chains and trees are stably realized in product sets of sufficient Newhouse thickness, and see \cite{McTBaire} for a topological variant, where the same authors show that all countably infinite bounded point configurations, including infinite trees, are stably realized in second category Baire sets  \cite{McTBaire} and link this area to the Erd\H{o}s similarity conjecture.

\vskip.125in 
 
In this paper, we begin by extending these types of results to generalized distance graphs, showing that 
arbitrary trees are stably realized in sets
$E$  with $\hd(E)$  sufficiently high (with threshold depending on $\phi$ but not on the tree $G$), 
%for an interval of values of $t$ of positive length, 
under the 
assumption that the generalized Radon transform associated with $\phi$ satisfies a suitable 
Sobolev mapping property.

\vskip.125in

We shall also show that our method allows one to prove such results for trees $G$ composed of elements of a fixed configuration;
for brevity, we illustrate this just for a tree of triangles.

\begin{example} Let  $\phi(u,v)=|u-v|$  denote the standard  Euclidean distance and
the tree $G$ be a path on three vertices, $\mathbb V=\{x^1,\, x^2,\, x^3\}$. 
In Section \ref{subsectionhypergraph} 
below, we show that the configurations consisting of $G$, or in fact any tree,  `decorated' with congruent triangles
can be stably realized in any
compact $E\subset {\Bbb R}^d$, $d \ge 4$, with $\hd(E)>(2d+3)/3$.
See Figure \ref{fig one} below,  the discussion in the next section, and further details in Section \ref{subsectionhypergraph}.
\end{example}

\vskip.125in 

\subsection{Structure of this paper} 
We begin by proving a result about tree structures in sets of a given Hausdorff dimension based on a 
general scheme that applies to a wide variety of situations. 
We then show that, if a (sufficiently symmetric)  configuration can be embedded in the distance graph
of a set $E$ and $\hd(E)$ is sufficiently high,  then  a {\it tree} of such configurations, 
\new{where the edges of the tree are  suitable subsets of the hyperedge defined by this 
configuration}, is also guaranteed to \new{be stably realizable} under a suitable dimensional threshold on $E$. 
We then provide concrete applications based on Fourier Integral Operator bounds for associated Radon transforms. 

\vskip.125in 

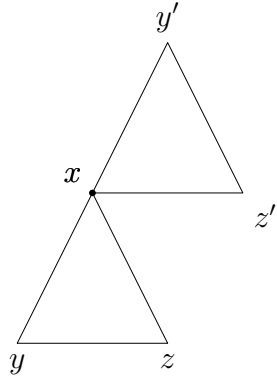
\begin{figure}
    \centering
    \begin{tikzpicture}
        % First triangle
        \tkzDefPoint(0,0){y}
        \tkzDefPoint(3.2,2.6){z}
        \tkzDefPoint(2,6){x}
        \tkzDrawPolygon(y,z,x)
        
        % Second triangle
       \tkzDefPoint(2,6){x}
        \tkzDefPoint(5.60,4.71){z'}
        \tkzDefPoint(8.05,8.03){y'}
        \tkzDrawPolygon(x,z',y')
        
        % Common vertex
        \tkzDrawPoints[fill=black](x,x)
         \tkzDrawPoints[fill=black](y,y)
          \tkzDrawPoints[fill=black](y',y')
        
        % Labels on vertices
        \tkzLabelPoint[below](y){$x^1$}
        \tkzLabelPoint[below](z){$y^1$}
        \tkzLabelPoint[above left](x){$x^2$}
       % \tkzLabelPoint[above left](x){$x$}
        \tkzLabelPoint[below right](z'){$y^2$}
        \tkzLabelPoint[above](y'){$x^3$}
        
          % Labels on edges
           \tkzLabelSegment[above left](y,x){$t^1$}
           \tkzLabelSegment[right](y,z){$t^3$}
            \tkzLabelSegment[right](x,z){$t^2$}
            
               \tkzLabelSegment[above left](x,y'){$t^1$}
           \tkzLabelSegment[below](x,z'){$t^3$}
            \tkzLabelSegment[right](y',z'){$t^2$}

        % Comment
      %  \node[align=center, below right] at (3.5, 1) {This picture shows a chain \\ of two triangles with a common vertex.};
    \end{tikzpicture}
    \caption{Chain of Two Congruent Triangles}
    \label{fig one}
\end{figure}

\subsection{Trees in general distance graphs} 
\label{treesubsection} 

\vskip.125in 

An important aspect of our approach is to formulate the realizability of trees (and certain other configurations) in a general setting. 
Let $\mu$ be a compactly supported nonnegative probability Borel measure 
on ${\Bbb R}^d$, and
$$K: {\Bbb R}^d \times {\Bbb R}^d \to [0,\infty)$$  a symmetric,
$(\mu\times\mu)$-integrable continuous function. 
Consider the graph, denoted by $G_K$, whose vertices are the points of $E \subset {\Bbb R}^d$, 
with two vertices $x,\, y\in E$ being connected by an edge iff $K(x,y)>0$. 

\vskip.125in 

The following result will allow us to reduce a variety of configuration problems to a series of concrete operator bounds.

\begin{theorem}[Tree building criteria] \label{maingeneral} 
Let $\mu$ and $K$ be as above. 
Define 
$$ U_Kf(x)=\int K(x,y) f(y) d\mu(y), \quad\forall f\in C_0(\mathbb R^d),$$ and suppose that 
\begin{equation} \label{assumption1}  \int \int K(x,y) d\mu(x)d\mu(y)>0, \end{equation} and 
\begin{equation} \label{assumption2}   \ U_K: L^2(\mu) \to L^2(\mu)\hbox{ continuously}. \end{equation} 

Then for any $k \ge 1$, 
\begin{equation} \label{chainlower} 
\int \cdots \int K(x^1,x^2)K(x^2,x^3) \cdots K(x^k,x^{k+1}) \,
d\mu(x^1)d\mu(x^2) \dots d\mu(x^{k+1})>0. \end{equation} 

\vskip.125in

More generally, let $T$ be a tree graph on $n$ vertices, $n \ge 2$,  with edge map ${\mathcal E}_T$. 
Define $K^*:\left(\R^d\right)^n\to [0,\infty)$ by
$$K^*(x^1,\dots,x^n)=\prod_{1\le i<j\le n,\, \mathcal E_T(i,j)=1} K(x^i,x^j).$$
Then
\begin{equation} \label{chainlowergeneral} 
\int \cdots \int K^*(x^1,\dots,x^n) \,
d\mu(x^1)d\mu(x^2) \dots d\mu(x^{n})>0. \end{equation} 

\end{theorem} 

In other words, the existence of trees even in this generalized setting is simply a consequence of assumptions (\ref{assumption1}) and (\ref{assumption2}). 
For example, in order to handle the Euclidean distance graph, where the vertices are points in a compact set $E$ and two vertices are connected by an edge 
iff the distance between them is equal to some fixed $t>0$, one takes $K=\sigma_t(x-y)$, 
where $\sigma$  is the  surface measure on the sphere of radius $t$. In this case, $U_K$ is the corresponding translation-invariant spherical averaging operator. 
As a technical point, even though $\sigma$ is a measure and not an $L^1$ function, 
the proof is accomplished by convolving $\sigma$ with the approximation to the identity at scale $\epsilon$ and checking that the estimates (upper and lower bounds (see e.g. \cite{BIT15}) do not depend on $\epsilon$. The existence of arbitrary trees in the case $K=\sigma_t(x-y)$ was previously established by the second and the third listed author in \cite{IT19}. 
Further, observe that non-degeneracy of the point configurations in question is guaranteed by the positivity of the integrals in  \eqref{chainlower} and \eqref{chainlowergeneral}, since degenerate configurations form  lower dimensional sets, 
which are of measure 0 with respect to $\otimes^n \mu$. 
\vskip.125in

The view point afforded by Theorem \ref{maingeneral}  also proves useful in the context of trees of hypergraphs, 
significantly expanding the scope of configurations that can be handled. 
For example, suppose that we want to show that a set $E$ contains many 2-chains of congruent triangles (see Fig. \ref{fig one}).
We are led to considering, for a Borel measure  $\mu$ supported on $E$, the expression
\begin{equation} \label{trianglehypergraph} \int^{(5)} F_1(x^1,x^2)F_2(x^2,y^1)F_3(y^1,x^1)F_1(x^2,x^3)F_2(x^3,y^2)F_3(y^2,x^2)\, 
d\mu(x^1)d\mu(y^1) d\mu(x^2) d\mu(y^2)d\mu(x^3), \end{equation} 
where, for $j=1,2,3$, the $F_j(x,y)$ are  non-negative $L^1$ functions which  will be smoothed out versions of 
$\delta\left(|x-y|-t_j\right)$, where the vector $\vec{t}=(t_1,t_2,t_3)$ can range over a set $S$ of side length vectors  
with nonempty interior in $\R^3_+$.
Then we may rewrite (\ref{trianglehypergraph}) in the form 
$$ \int \int \int K(x^1,x^2)K(x^2,x^3)\,  d\mu(x^1)d\mu(x^2)d\mu(x^3),
%= \int \int \int K(y,x)K(x,y') d\mu(y)d\mu(x)d\mu(y'),
$$ 
where
\begin{equation} \label{trianglekernel} K(x,y):=F_1(x,y) \int F_2(x,z)F_3(y,z) d\mu(z)=K(y,x) \end{equation} 
is symmetric and satisfies (\ref{assumption1}) and (\ref{assumption2}).
Thus, Theorem \ref{maingeneral}  applies,
establishing the existence in $E$  of chains of two congruent triangles, 
for all vectors $\vec{t}\in S$.
More details can be found in Section \ref{subsectionhypergraph} below,
together with the extension from 2-chains to arbitrary trees of triangles and certain other configurations.

\vskip.125in 

\section{Consequences of Theorem \ref{maingeneral}} 

We begin with corollaries in the setting of trees, followed by hypergraphs, showing that 
Theorem \ref{maingeneral} reduces the existence of a wide variety of configuration to the verification of conditions (\ref{assumption1}) 
and (\ref{assumption2}). These conditions amount to certain function space estimates which may, depending on the 
particular result, be of 
greater or lesser difficulty to establish. 
We will now proceed to work out a variety of such examples. 

\vskip.125in 

\subsection{Generalized Radon transforms} 

\begin{corollary}[Realizing trees in sets of sufficient Hausdorff dimension] \label{main} 
Let $T$ be a tree graph on $n$ vertices, and  ${\mathcal E}_T$ its edge map. 
Let $\phi:{\Bbb R}^d \times {\Bbb R}^d \to {\Bbb R}$ be  continuous. 
Suppose that for all $t>0$,  $\phi$ is smooth near $\Sigma_t:=\{(x,y): \phi(x,y)=t\}$, with $\nabla_x \phi(x,y), \nabla_y \phi(x,y) \not=\vec{0}$, 
so that $\Sigma_t\subset \mathbb R^d\times \mathbb R^d$ is  smooth and for each $x\in\mathbb R^d$, $\Sigma_t^x:=\{y \in {\Bbb R}^d: \phi(x,y)=t \}\subset\mathbb R^d$ is smooth.
Further assume that,
if $\psi$ is a smooth cut-off and $\sigma_{x,t}$ is the surface measure on 
$\Sigma_t^x$,  there is some $\alpha>0$ such that the generalized Radon transform
$$ R^{\phi}_tf(x):=\int_{\Sigma_t^x} f(y) \, \psi(y)\,  d\sigma_{x,t}(y),$$ 
is continuous $L^2({\Bbb R}^d) \to L^2_{\alpha}({\Bbb R}^d)$, locally uniformly in $t$.

Then, if $E \subset {\Bbb R}^d$ is compact with Hausdorff dimension $\hd(E)>d-\alpha$,
$T$ is stably realizable in $E$.
I.e., there is a non-empty open interval $I\subset\R_+$ such that
$T$ is realizable in  $E$ with any gap $t\in I$: 
for all $t\in I$ there exist distinct $x^1, \dots, x^{n} \in E$ such that 
$\phi(x^i,x^j)=t$ for $(i,j)$ such that ${\mathcal E}_T(i,j)=1$. \end{corollary}

\vskip.125in

(We 
refer to \cite{T11,Graf14} for treatments of the $L^2$-based Sobolev spaces, $L^2_\alpha,\, \alpha\in\R$, 
that we will use.)
Corollary \ref{main} will be proved in Section \ref{sec proof of cor}.
As an example, if $\phi(x,y)=|x-y|$ is the Euclidean distance,
the corresponding $R^{\phi}_t$ are spherical means operators which are smoothing of order $\alpha=(d-1)/2$, and thus
 the conclusions of Cor. \ref{main} hold if the Hausdorff dimension of $E$ is greater than $\frac{d+1}{2}$. 
 More generally, this $L^2$ Sobolev regularity holds within the conjugate locus for the distance induced by any Riemannian metric, leading immediately to the following.

\vskip.125in

\begin{corollary}[Realizing trees on Riemannian manifolds]  \label{riemann} Let $(M,g)$ be a compact Riemannian manifold without  boundary or 
conjugate points, of dimension $d \ge 2$. 
Let $E$ be a subset of $M$ of Hausdorff dimension $>\frac{d+1}{2}$, and let $\rho_M$ denote the 
induced Riemannian distance function on $M$.

Let $T$ be any tree graph on $n$ vertices, and let ${\mathcal E}_T$ be the corresponding edge map.
Then $T$ is realizable in $E$, in the sense that there exist $x^1, \dots, x^{n} \in E$ and a non-empty interval $I$ such 
that $\rho_M(x^i,x^j)=t$ for $(i,j)$, $t \in I$, such that ${\mathcal E}_T(i,j)=1$. 
\end{corollary}

\vskip.125in

\subsection{Trees of triangles} \label{subsectionhypergraph} 
We now illustrate  applications of Theorem \ref{maingeneral}
to more complicated configurations,
based on the approach briefly
described below the statement of Theorem  \ref{maingeneral}
for 2-chains of congruent triangles, 
and this is where we pick up the narrative and provide more details
for some specific examples. 
\smallskip

\new{Theorem \ref{maingeneral}  yields the existence of trees of $k$-point configurations
for a wide variety of configurations studied in the literature on what are now called Mattila-Sj\"olin functions.
For $d,\, p\in\mathbb N$ and $k\ge 2$, let $\Phi:\left(\R^d\right)^k \to \R^p$ be a 
continuous function which is smooth (except possibly on a lower dimensional set).
For a compact $E\subset\R^d$, the {\it $\Phi$-configuration set} of $E$ \cite{GGIP12} is the compact set
$$\Delta_\Phi(E):=\left\{\, \Phi\left(x^1,\dots,x^k\right): \, x^1,\dots,x^k\in E\, \right\}\subset \R^p.$$
Then $\Phi$ is said to be a {\it Mattila-Sj\"olin function} if there is some $s_{\Phi}<d$ such that $\hd(E)>s_\Phi$ ensures that
$\Delta_\Phi(E)$ has nonempty interior.
See \cite{MS99,IMT11,GIT19,GIT20,PRA21,KPS21,GIT22,PRA22,GGPP23} for results of this type.} 

Our approach in \cite{GIT19,GIT20,GIT22} was as follows:
Let $\mu$ be a Frostman measure supported on $E$ and of finite $s$-energy.
Then the Radon-Nikodym derivative of the
{\it configuration measure}
$\nu_\Phi:=\Phi_*\left(\mu\otimes\cdots\otimes\mu\right)$ on $\R^p$ with respect to Lebesgue  measure $d\bf t$
can be represented as a multilinear form $\Lambda_{\mathbf t}(\mu,\dots,\mu)$ ,
with multilinear kernel
$$L_{\mathbf t}(x^1,\dots,x^k)=\delta\left(\Phi\left( x^1,\dots,x^k\right)-{\mathbf t}\right),$$
where $\delta$ is the Dirac delta at $0\in \R^p$. 
The method of partition optimization \cite{GIT19,GIT20} and its local and microlocal variants \cite{GIT22} 
allow one to obtain estimates of the form
$$\left| \Lambda_{\mathbf t}\left(f_1,\dots, f_k\right)\right| \le C \prod_{j=1}^k || f_j ||_{L^2_{r_j}}$$
for negative $r_j$ with a lower bound on $\sum r_j$ depending on $\Phi$.
Applying this to $f_1=\cdots=f_k=\mu$ having finite $s$-energy (which implies that $\mu \in L^2_{\frac{s-d}2}$)
yields that, for $\mathbf t$ in a set $S\subset\R^p$ with nonempty interior,
\begin{equation}\label{eqn Lambda}
0<\Lambda_{\mathbf t}(\mu,\mu,\dots,\mu):=\langle\, \delta\left(\Phi\left( x^1,\dots,x^k\right)-{\mathbf t}\right),\, \mu\otimes\cdots\otimes\mu\, \rangle = C<\infty
\end{equation}
if $\hd(E)>s_\Phi$, and furthermore $\Lambda_{\mathbf t}(\mu,\dots,\mu)$ is continuous in ${\mathbf t}\in S$.
Here, the pairing $\langle \cdot,\, \cdot \rangle$ is between distributions and Sobolev functions on $\R^{kd}$. 

The left side of \eqref{eqn Lambda} can be rewritten in various equivalent ways by partitioning the variables and integrating out some first.
In particular, fixing distinct indices $1\le i<j\le k$, for  points $x,y\in \R^d$,  let  
$$\hat{\mathbf x}^{ij}=\left(x^1,\dots,x^{i-1},x,x^{i+1},\dots,x^{j-1},y,x^{j+1},\dots x^k\right).$$
Then \eqref{eqn Lambda} implies  that
\be\label{eqn partition}
K(x,y):=\big\langle\, \delta\left( \Phi\left(\hat{\mathbf x}^{ij}\right) \right),\, 
\bigotimes_{\buildrel{l=1}\over{ l\ne i,j}}^k \, \mu\left(x^l\right)\, \big\rangle
\ee
has the property that 
\begin{equation}\label{eqn 1.3 again}
0<\int\int K(x,y)\mu(x)\mu(y)=  C<\infty,
\end{equation} which is  condition \eqref{assumption1} from Theorem \ref{maingeneral}.
The following argument shows that, by replacing $\mu$ with its restriction to a chosen subset $F\subset E,\, \mu(F)>0$, one can preserve
\eqref{assumption1} while also ensuring that \eqref{assumption2} holds.
\medskip

To start, from \eqref{eqn 1.3 again} it follows that
$$\mu\left\{y\,:\, \int K(x,y)\, d\mu(x) >3C\right\} \le \frac13.$$
Thus, if we define 
$$F_1=\left\{y\in E\, :\, \int K(x,y)\, d\mu(x)\le 3C\, \right\},$$
then $\mu(F_1)\ge 2/3$. Denoting by $\mu_1$ the restriction of $\mu$ to $F_1$,
it follows that 
\begin{equation}\label{eqn schur1a}
\int K(x,y)\, d\mu(x) \le 3C\hbox{ for all } y\in F_1,
\end{equation}
and
$$\int\int K(x,y)\, d\mu(x)\, d\mu_1(y) \le C.$$
Using the last inequality and the positivity of the integrand, we can change the order of integration and 
repeat this argument with respect to $x$. Since
$$\mu\left\{x\,:\, \int K(x,y)\, d\mu_1(y) >3C\right\} \le \frac13,$$
the set $F_2: =  \left\{x\,:\, \int K(x,y)\, d\mu_1(y) \le3C\right\}$ has $\mu(F_2)\ge 2/3$.
Denoting $\mu|_{F_2}$ by $\mu_2$, we then have
\begin{equation}\label{eqn schur2a}
\int K(x,y)\, d\mu_1(y) \le 3C\hbox{ for all } x\in F_2.
\end{equation}
$$\mu(F_1\cap F_2)\ge 1-\frac13-\frac13=\frac13,$$
and we denote $F_1\cap F_2$ by $\tilde{E}$ and $\mu|_{\tilde{E}}$ by $\tilde\mu$.
Finally, we have
\begin{equation}\label{eqn schur}
\int K(x,y)\, d\tilde\mu(x) \le 3C\hbox{ for all } y\in \tilde{E},\quad \int K(x,y)\, d\tilde\mu(y) \le 3C\hbox{ for all } x\in \tilde{E},
\end{equation}
so that Young's inequality applies to the integral kernel $\tilde{K}=K|_{ \tilde{E}\times  \tilde{E}}$,
which thus defines a bounded operator $U_{\tilde{K}}:L^2(\tilde\mu)\to L^2(\tilde\mu)$.
Erasing all the tildes,  we see that  \eqref{assumption2} is satisfied.
\medskip

Finally, in order to apply Thm. \ref{maingeneral}  to obtain 
trees of $\Phi$-configurations, 
we need that that $K(x,y)$ be symmetric in $x$ and $y$, 
and for this one needs to impose some symmetry conditions on the configuration function $\Phi$.
For simplicity, take $i=1,\, j=2$ in \eqref{eqn partition}; then we demand that for some $A\in GL(p,\R)$ and some permutation
$\mathbf{\tilde{x}}^{12}$ of the $k-2$ variables in $\mathbf{\hat{x}}^{12}$, we have
\be\label{eqn equiv partition}
\Phi\left(y,x,\mathbf{\hat{x}}^{12}\right)=A\circ \Phi\left(x,y,\mathbf{\tilde{x}}^{12}\right),
\ee
so that $K(y,x)=c\cdot K(x,y)$, with $c=|A|^{-1}$, which preserves \eqref{eqn 1.3 again} and  is good enough for our purposes.

\bigskip

\subsubsection{Building a tree of congruent triangles}
The mechanism of this  paper applies to any configuration for which we can prove that the natural measure 
associated to  the 
configuration satisfies the assumptions of Theorem 
\ref{maingeneral}, in the sense described  below the statement of that result. 
We give just one illustrative example, namely the existence of trees of congruent  triangles in $E$,
using the following result from \cite{PRA21}; see also \cite{GIT22} for an alternate proof using microlocal analysis. 

\begin{theorem} \label{microlocalpaper} (\cite{PRA21}) If $E \subset \mathbb{R}^d$, $d \geq 4$, is compact with $\text{dim}_H(E) > \frac{2d + 3}{3}$, 
then the set of congruence classes of triangles with vertices in $E$,
\begin{equation}
    \left\{(|x-y|, |x-z|, |y-z|) : x, y, z \in E\right\},
\end{equation} has nonempty interior in $\mathbb{R}^3$.
\end{theorem}

Moreover, it is shown that the natural measure supported on 
$$   \left\{(|x-y|, |x-z|, |y-z|) : x, y, z \in E \right\},$$ 
namely the configuration measure $\nu_{triangle}$ defined by 
$$ \int f(t^1,t^2,t^3) d\nu_{triangle}(t^1,t^2,t^3):=\int \int \int f(|x-y|, |x-z|, |y-z|) d\mu(x) d\mu(y) d\mu(z), $$ is continuous away from the degenerate 
triangles,  and from this the conclusion of Theorem \ref{microlocalpaper} is ultimately obtained. 

 Given  a side length vector  $\vec{t}=(t^1,t^2,t^3)$ in the non-empty interior, say $S$,  of the configuration set guaranteed 
 by Theorem \ref{microlocalpaper}, we can build a tree of congruent
 triangles with side lengths $\vec{t}$, with any two triangles joined at exactly one vertex, as follows. 
At such an $\vec{t}$ the measure $\nu_{triangle}$ has a continuous density function, namely 
$$ \nu_{triangle}= \lim_{\epsilon \to 0^{+}} \nu_{triangle}^{\epsilon},$$ where 
$$ \nu_{triangle}^{\epsilon}\left(\, \vec{t}\, \right)=  \int \int \int \sigma_{t^1}^{\epsilon}(x-y) \sigma^{\epsilon}_{t^2}(x-z) \sigma^{\epsilon}_{t^3}(y-z) 
d\mu(x) d\mu(y) d\mu(z).$$ 

We now define the approximate kernel 
$$ K_{\vec{t}}^{\epsilon}(x,y)=\sigma_{t^1}^{\epsilon}(x-y) \int \sigma^{\epsilon}_{t^2}(x-z) \sigma^{\epsilon}_{t^3}(y-z) d\mu(z),$$ 
which satisfies the equivariance property \eqref{eqn equiv partition} with $A\in GL(3,\R)$ interchanging $t^2$ and $t^3$.
This is the object to which we apply Theorem \ref{maingeneral}, and the argument is complete. 
Following the same procedure we can 
produce an arbitrary tree of triangles, not just two triangles joined at a vertex. 
\medskip

\subsubsection{Trees of equi-area triangles} This method can be applied to obtain  the existence of arbitrary trees 
for some, but not all, of the $k$-point configurations for which nonempty interior of configurations sets were established
 in \cite{GIT20,GIT22}.  One of these concerned areas of triangles in the plane:
 
 \begin{theorem}\label{thm areas} ( \cite[Thm. 1.1(i)]{GIT20})
If $E\subset\R^2$ is compact with $\hd(E)>5/3$, then 
 the set  of signed areas of triangles determined by triples of points of $E$, 
\be\label{def areas}
\left\{\frac12\det\left[x-z,\, y-z\right]\, : x,y,z\in E\right\}\subset \R,
\ee
contains an open interval.
\end{theorem}

The $\R^1$-valued configuration function $\Phi$ of three variables in $\R^2$, 
$\Phi(x,y,z)=\det\left[x-z,\, y-z\right]$ satisfies \eqref{eqn equiv partition} with factor $-1$,
so that the method above applies. Hence, for $\hd(E)>5/3$ and  for an arbitrary tree $T$, and areas $A$ in an open interval,
there exist copies of $T$ in $E$ and auxiliary points $y^{ij}\in E$ for each
$(i,j)$ with ${\mathcal E}_T(i,j)=1$, such that $x^i,\, x^j$ and $y^{ij}$ span a triangle of area $A$.

\section{Proof of Theorem \ref{maingeneral}} 

\vskip.125in 

Our basic scheme is the following. We first prove the result for paths with $k=2^m$ vertices by utilizing Cauchy-Schwarz and the assumption 
(\ref{assumption1}). We then induct downwards to fill in the gaps between the dyadic numbers after first pigeonholing to a subset where $U_K1$ is not too 
large. Finally, we notice that the flexibility afforded by our arguments allows us to extend the case of a path to a general tree. 

\subsection{The case $k=2^m$} 
\label{subsectiondyadic}

We begin by proving (\ref{chainlower}). Set 
$$c_{lower}:=\int \int K(x,y) d\mu(x)d\mu(y),$$
which is $>0$ by \eqref{assumption1},
and define
$$ C_k(\mu)=\int \dots \int \prod_{j=1}^k K(x^j,x^{j+1}) \prod_{i=1}^{k+1} d\mu(x^i).$$

Suppose that $k=2^m$ for some $m$. Then, by repeated application of Cauchy-Schwarz, we have
$$ C_{2^m}(\mu)=\int \dots \int \prod_{j=1}^k K(x^j,x^{j+1}) \prod_{i=1}^{k+1} d\mu(x^i)$$
$$ \ge {\left( \int \int K(x,y) d\mu(x)d\mu(y) \right)}^{2^m}>c_{lower}^{2^m}>0,$$
where we used (\ref{assumption1}) and the assumption that $\mu$ is a probability measure.

\vskip.125in 

\subsection{Refinement to a subset where $U_K1$ is not too large} 
\label{subsectionrefinement} 

In order to deal with general $k$, we need to do a bit of pigeon-holing. Observe that,
if $C_{norm}:=\left|\left| U_K\right|\right|_{L^2\to L^2}$,  
$$ \mu \{x: \left(U_K1\right)(x)>\lambda \} \leq \frac{1}{\lambda^2} \int {|U_K1(x)|}^2 d\mu(x) \leq \frac{C_{norm}^2}{\lambda^2}$$ by (\ref{assumption2}). 
It follows that if $\lambda=NC_{norm}$, with $N>2$ to be determined later, then 
\begin{equation} \label{upperboundpigeonhole} \left(U_K1\right)(x) \leq NC_{norm} \ \text{on a set} \  E' \  \text{with} \ \mu(E') \leq \frac{1}{N^2}. \end{equation} 

If we replace the constant function $1$ in (\ref{upperboundpigeonhole}) by the indicator function of $E'$, the upper bound still holds. Moreover, if we let $\mu'$ denote $\mu$ restricted to $E'$, we have 
$$ \int \int K(x,y) d\mu'(x) d\mu'(y)$$
$$=\int \int K(x,y) d\mu(x) d\mu(y)- \int \int K(x,y) f(x)f(y) d\mu(x)d\mu(y)$$
$$=I-II,$$ where $f$ is the indicator function of the set where $U_K1(x)>2C_{norm}$. By assumption, $I=c_{lower}>0$. Observe that 
$$ II \leq C_{norm} \cdot  {||f||}^2_{L^2(\mu)} \leq \frac{C_{norm}^3}{\lambda^2} \leq \frac{C_{norm}}{N^2} \leq \frac{c_{lower}}{2}$$ 
if we choose
$$ N \ge \sqrt{\frac{2C_{norm}}{c_{lower}}}.$$ 

With a slight abuse of notation, we can now rename $\mu'$ back to $\mu$, renormalize, and pretend that from the very beginning we had a set $E$, equipped with the Borel measure $\mu$, such that $\left(U_K1\right)(x)$ is bounded above by some uniform constant $C$, and both (\ref{assumption1}) and (\ref{assumption2}) hold. 

\vskip.125in 

\subsection{Paths of arbitrary finite length and transition to trees} Using the results just obtained in Section \ref{subsectionrefinement}, 
one sees that for $k \ge 2$, 
$$ C_{k}(\mu) \leq C \cdot C_{k-1}(\mu),$$ where $C$ is the upper bound on $\left(U_K1\right)(x)$.

Proceeding by induction we get a lower bound on a path of arbitrary length. In particular, and this notion will come in handy in a moment, 
having built a path in $E$  with $2^m$ links, we have also built a path of smaller length. 

In order to build an arbitrary tree, we use the simple principle that if $T,T'$ are trees, and $T$ is contained in $T'$, then building $T'$ in $E$ 
implies that we can build $T$. Given a tree $T$, let 
\begin{equation} \label{treeexpression} T(\mu)=\int \dots \int \prod_{(i,j) \in {\mathcal E}_T} K(x^i,x^{j}) \prod_{i=1}^{k+1} d\mu(x^i), \end{equation}
where ${\mathcal E}_T$ is as above.  

\vskip.125in 

We shall need the following definition. 

\begin{definition}[A wrist of a tree] \label{wristdef}  Let G be a connected tree graph. We say that $w$ is a wrist of order $n$ if the following conditions hold: 

i) $w \in V$, the vertex set of $G$. 

ii) $V=V_1 \cup V_2$, where $V_1 \cap V_2=\{w\}$.  

iii) Vertices from $V_1 \backslash \{w\}$ are not connected by edges to vertices in $V_2 \backslash \{w\}$. 

iv) Let $G_1$ denote $G$ restricted to $V_1$. Then $G_1$ is the union of finitely many chains 

$C_1, C_2, \dots, C_n$ such that vertices of $C_i$ only intersect vertices of $C_j$, $i \not=j$, at $w$.
\end{definition}

\begin{example} i) Consider a chain on three vertices, with vertices $v_1, v_2, v_3$ such that $v_1$ is connected by an edge to $v_2$, and $v_2$ is 
connected by an edge to $v_3$, but $v_1$ and $v_3$ are not connected. Then $v_1$ is a wrist because there is a chain with vertices $v_1, v_2, v_3$ with 
one endpoint at $v_1$. The vertex $v_3$ is a wrist for the same reason. The vertex $v_2$ is also a wrist because two chains, namely the one with vertices 
$v_2, v_1$, and the one with vertices $v_2, v_3$ have $v_2$ as an endpoint. 

\vskip.125in 

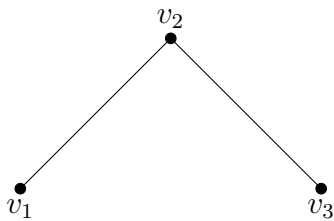
\begin{tikzpicture}
    % Define vertices
    \coordinate (v1) at (0,0);
    \coordinate (v2) at (2,2);
    \coordinate (v3) at (4,0);
    
    % Draw vertices
    \filldraw (v1) circle (2pt) node[below] {$v_1$};
    \filldraw (v2) circle (2pt) node[above] {$v_2$};
    \filldraw (v3) circle (2pt) node[below] {$v_3$};
    
    % Draw segments
    \draw (v1) -- (v2);
    \draw (v2) -- (v3);
\end{tikzpicture}

\vskip.125in 

ii) Consider a complete graph on three vertices. Then no vertex is a wrist. 

\vskip.125in 

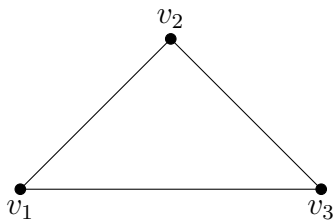
\begin{tikzpicture}
    % Define vertices
    \coordinate (v1) at (0,0);
    \coordinate (v2) at (2,2);
    \coordinate (v3) at (4,0);
    
    % Draw vertices
    \filldraw (v1) circle (2pt) node[below] {$v_1$};
    \filldraw (v2) circle (2pt) node[above] {$v_2$};
    \filldraw (v3) circle (2pt) node[below] {$v_3$};
    
    % Draw segments
    \draw (v1) -- (v2);
    \draw (v2) -- (v3);
    \draw (v1) -- (v3); % Add this line to connect v1 and v3
\end{tikzpicture}

\vskip.125in 

iii) Consider a graph on four vertices $v_1, v_2, v_3, v_4$, where $v_1$ and $v_2$ are connected, $v_2$ and $v_3$ are connected, $v_2$ and $v_4$ are 
connected, and there are no other edges. Then $v_2$ is the only wrist. 

\vskip.125in 

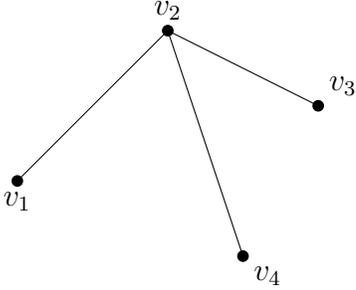
\begin{tikzpicture}
    % Define vertices
    \coordinate (v1) at (0,0);
    \coordinate (v2) at (2,2);
    \coordinate (v3) at (4,1);
    \coordinate (v4) at (3,-1);
    
    % Draw vertices
    \filldraw (v1) circle (2pt) node[below] {$v_1$};
    \filldraw (v2) circle (2pt) node[above] {$v_2$};
    \filldraw (v3) circle (2pt) node[above right] {$v_3$};
    \filldraw (v4) circle (2pt) node[below right] {$v_4$};
    
    % Draw segments
    \draw (v1) -- (v2);
    \draw (v2) -- (v3);
    \draw (v2) -- (v4);
\end{tikzpicture}

\vskip.125in 

\end{example}

Our argument is based on the fact that every tree which is not a chain contains a wrist of order $>1$. 

\begin{lemma}[Any nontrivial tree contains a wrist] \label{wristsexistlemma} Let $T$ be a finite connected tree. Then either $T$ is a chain, or $T$ contains a wrist of order $>1$. \end{lemma} 

To prove the lemma, let $v_1, v_2, \dots, v_m$ denote the (distinct) vertices of degree $1$ in $T$. We move from each $v_j$ until we encounter a vertex of degree $\ge 3$. If such a vertex does not exist, then $T$ is clearly a chain. In this way, we assign a vertex $w_j$, of degree $\ge 3$, to each $v_j$. We claim that there exist $i,j$, $i \not=j$, such that $w_i=w_j$. Suppose not. Remove all the $v_j$s and the vertices and edges that lead up to, but not including, $w_j$. The resulting graph $T'$ is still a connected tree. Each vertex $w_j$ in $T'$ has a degree $\ge 3-1=2$, so the vertices of degree $1$ in $T'$ are not any of the $v_j$s or any of the $w_j$s. This means that those vertices were present in the original tree graph $T$, but this is impossible since we removed them. 

Now that we have shown that there exists $i \not=j$ such that $w_i=w_j$, it is not difficult to see that this $w_i$ is a wrist of order $>1$, as desired. This completes the proof of the lemma.

\vskip.125in 

Let $w_0$ denote a wrist point in $T$, which, as we just proved, is guaranteed to exist. We now rewrite (\ref{treeexpression}) in the form 

\begin{equation} \label{preholder} \int U_K(w) C_1(w) C_2(w) \dots C_k(w) dw, \end{equation} where 

$$ C_j(w)=\int \dots \int K(w,x^{j,1})K(x^{1,1},x^{j,2}) \dots K(x^{j,n_j-1},x^{j, n_j}) dx^{j,1} \dots dx^{j,n_j}.$$

Adding vertices and edges, if necessary, we can make all the chains have the same length, $n_{max}$. 
One can then estimate (\ref{treeexpression}) using H\"{o}lder's inequality. 

$$ \int U_K(w) {\left( \int \dots \int K(w,x^1) K(x^1,x^2) \dots K(x^{n_{max}-1}, x^{n_{max}}) dx^1 \dots dx^{n_{max}} \right)}^k dw$$ 
$$ =\int {\left( \int \dots \int K(w,x^1) K(x^1,x^2) \dots K(x^{n_{max}-1}, x^{n_{max}}) dx^1 \dots dx^{n_{max}} \right)}^k U_K(w)dw$$ 
$$ \ge \frac{{\left( \int  \int \dots \int K(w,x^1) K(x^1,x^2) \dots K(x^{n_{max}-1}, x^{n_{max}}) dx^1 \dots dx^{n_{max}} U_K(w)dw \right)}^k}{{\left( \int U_K(w) dw\right)}^{k-1}}$$
\begin{equation} \label{treereduced} \ge C {\left( \int  \int \dots \int K(w,x^1) K(x^1,x^2) \dots K(x^{n_{max}-1}, x^{n_{max}}) dx^1 \dots dx^{n_{max}} U_K(w) dw \right)}^k \end{equation} since we have an upper bound for $\int U_K(w) dw$ by a repeated use of (\ref{upperboundpigeonhole}). 

\vskip.125in 

In other words, 
$$ T(\mu) \ge c \cdot T'(\mu),$$ 
where $T'$ is the tree obtained from $T$ by removing all but the longest chain emanating from the wrist $w$. 
It is clear that $T'$ has fewer vertices (and hence edges) than $T$. 
Proceeding in this way shows that given any tree $T$, 
there exists a tree $T^{*}$ containing $T$ and a positive constant $c^{*}$ such that 
$$ T^{*}(\mu) \ge c^{*} \cdot c_{lower}>0.$$ 

This completes the proof of Theorem \ref{maingeneral}. 

\vskip.125in

\section{Proof of Corollary \ref{main}}\label{sec proof of cor}

\vskip.125in

The proof of Corollary \ref{main} follows, in view of Theorem \ref{maingeneral}, from the following results. The first one follows from the proof of the main result in \cite{GIT19} (also see \cite{GIT20}).

\begin{theorem}[Establishing the lower bound \eqref{assumption1}]  \label{lowerboundtheorem} 
Let $\phi$, $R_t^{\phi}$ be as in the statement of Theorem \ref{main}, with
$R^{\phi}_t: L^2({\Bbb R}^d) \to L^2_{\alpha}({\Bbb R}^d)$  for some $\alpha>0$,
uniformly for $t$ in a non-trivial interval  $I_0\subset\R$.
Let $E$ be
a compact set of Hausdorff dimension $\hd(E)>d-\alpha$,
and $\mu$  a Frostman measure on $E$ of finite $s$-energy for some $s>d-\alpha$.
Then
$$J(t):= \int R_t^{\phi}\mu(x) d\mu(x)$$ is a continuous function on $I_0$, and there exists a non-empty 
 open  interval $I\subseteq I_0$ and a $c_\delta>0$ such that {for all $t\in I$,}
\begin{equation} \label{lowerbound} \int R_t^{\phi}\mu(x) d\mu(x) \ge c_{\delta}>0. \end{equation}
\end{theorem} 

Theorem \ref{lowerboundtheorem} establishes that assumption \eqref{assumption1} in 
Theorem \ref{maingeneral} is satisfied (with $K(x,y)$ the Schwartz kernel of $R_t^\phi$) uniformly for $t\in I$. 

\begin{theorem}[Establishing the upper bound \eqref{assumption2}] \label{L2theorem} Let $\phi$, $R_t^{\phi}$ be as in the statement of Corollary \ref{main}. 
Suppose that for some $\alpha>0$ and all $t>0$, $R^{\phi}_t: L^2({\Bbb R}^d) \to L^2_{\alpha}({\Bbb R}^d)$ is bounded. 
Let $E\subset\mathbb R^d$ be compact, with Hausdorff dimension  greater than 
$d-\alpha$, and $\mu$  a Frostman measure on $E$. Then for any $t>0$,
\begin{equation} \label{L2upperboundest} {||R_t^{\phi}f||}_{L^2(\mu)} \leq K {||f||}_{L^2(\mu)}. \end{equation} 
\end{theorem}

This establishes the assumption (\ref{assumption2}) in Theorem \ref{maingeneral}. 

\vskip.125in

\begin{remark} The constant $K$ above only depends (uniformly) on the implicit constants in the assumptions of Corollary \ref{main}. By standard FIO theory (see e.g. \cite{GS94}, Section 2, for similar calculations) $K$ depends only on the ambient dimension $d$, the Hausdorff dimension of the support of $\mu$, the bounds implicit in the Sobolev estimate for $R_t^{\phi}$ and the Frostman constant, i.e., the constant such that $\mu(B(x,r)) \leq Cr^{\alpha}$ for any $\alpha<d$, where $B(x,r)$ is the ball of radius $r$ (sufficiently small) centered at $x$ in the support of $\mu$. \end{remark}

\vskip.125in

\subsection{Proof of Theorem \ref{lowerboundtheorem}}\label{lowerboundtheorem}

This is essentially proven in \cite{GIT19}, but  for the sake of completeness we include the argument.
The assumption  $\nabla_x \phi(x,y), \nabla_y \phi(x,y) \not=\vec{0}$ on $\{(x,y): \phi(x,y)=t\}$
allows one to conjugate by an elliptic pseudodifferential operator of any order $r\in\R$, so that
$R^{\phi}_t: L^2_r({\Bbb R}^d) \to L^2_{r+\alpha}({\Bbb R}^d)$, locally uniformly in $t$. 
Since $\mu$ has finite $s$-energy, $\mu\in L^2_{{s-d}/2}$. 
Thus $R^{\phi}_t\mu\in L^2_{(s-d)/2+\alpha}$, 
and this will pair boundedly against $\mu$ if $(s-d)/2+\alpha+(s-d)/2\ge 0$, i.e., if 
$s>d-\alpha$. Furthermore, by continuity of the integral, this is continuous in $t$.
Since the integral of $J$ in $t$ is positive by the coarea formula, there must be a $t_0$ at which $J(t_0)>0$,
and hence there is a non-empty open interval on which $J$ is strictly positive.

\vskip.125in

\subsection{Proof of Theorem \ref{L2theorem}}
\label{subsectionL2theorem}

It is enough to show that
$$ \int R_t^{\phi}f\mu(x) g(x) d\mu(x) \leq C<\infty$$ for $g$ such that ${||g||}_{L^2(\mu)}=1$. Let ${(f\mu)}_j$ denote the Littlewood-Paley piece of $f\mu$ on scale $j \ge 0$. Negative scales are straightforward and will be handled separately. We are going to bound
$$ \langle R_t^{\phi}{(f\mu)}_j, {(g \mu)}_{j'}\rangle\, = \,\langle\widehat{R_t^{\phi}{(f\mu)}_j}, \widehat{{(g\mu)}_{j'}}\rangle.$$

By a standard orthogonality argument for generalized Radon transforms, the expression above decays rapidly when $|j-j'| \ge 5$. It follows that it suffices to bound
$$ \sum_{|j-j'| \ge 5} <\widehat{R_t^{\phi}{(f\mu)}_j}, \widehat{{(g\mu)}_{j'}}>$$
$$ \leq \sum_{|j-j'| \ge 5} {\left( \int {|\widehat{R_t^{\phi}{(f\mu)}_j}(\xi)|}^2 d\xi \right)}^{\frac{1}{2}} \cdot
{\left( \int {|\widehat{ {(g\mu)}_{j'}(\xi)}|}^2 d\xi \right)}^{\frac{1}{2}}$$
\begin{equation} \label{almost} \leq C \sum_{|j-j'| \ge 5} 2^{-j\alpha} {\left( \int {|\widehat{ {(f\mu)}_{j}(\xi)}|}^2 d\xi \right)}^{\frac{1}{2}} {\left( \int {|\widehat{ {(g\mu)}_{j'}(\xi)}|}^2 d\xi \right)}^{\frac{1}{2}}. \end{equation}

We shall need the following basic estimate. See \cite{IKSTU19} and \cite{EIT11} for similar results.

\begin{lemma} \label{schurtest} With the notation above, for any $\epsilon>0$,
$$ {\left( \int {|\widehat{ {(f\mu)}_{j}(\xi)}|}^2 d\xi \right)}^{\frac{1}{2}} \leq C_{\epsilon} 2^{j(d-s+\epsilon)/2} {||f||}_{L^2(\mu)}.$$
\end{lemma}

With Lemma \ref{schurtest} in tow, the expression in (\ref{almost}) is bounded by
$$ C 2^{-j \alpha} 2^{j(d-s+\epsilon)} {||f||}_{L^2(\mu)} {||g||}_{L^2(\mu)},$$ so the sum over $j$ is bounded by $C{||f||}_{L^2(\mu)}$ provided that $s>d-\alpha$, as claimed. This completes the proof of Theorem \ref{L2theorem}, once we establish Lemma \ref{schurtest}.

To prove Lemma \ref{schurtest}, we write
$$  \int {|\widehat{ {(f\mu)}_{j}(\xi)}|}^2 d\xi= \int {|\widehat{f\mu}(\xi)|}^2 \psi^2(2^{-j}\xi) d\xi,$$ 
where $\psi$ is a smooth cut-off function supported in the annulus
$$\left\{\xi \in {\Bbb R}^d: \frac{1}{2} \leq |\xi| \leq 4 \right\}.$$

This expression is bounded by
$$  \int {|\widehat{f\mu}(\xi)|}^2 \widehat{\rho}(2^{-j} \xi) d\xi,$$ where $\rho$ is a suitable cut-off function.

By Fourier inversion and a limiting argument (see \cite{W04}), this expression equals
$$ 2^{dj} \int \int \rho(2^j(x-y)) f(x)f(y) d\mu(x) d\mu(y)=<U_jf,f>,$$ where
$$ U_jf(x)=\int 2^{dj} \rho(2^j(x-y)) f(y) d\mu(y),$$ and $\langle \, ,\, \rangle$ is the $L^2(\mu)$ inner product.

Since
$$ \int 2^{dj} \rho(2^j(x-y)) d\mu(y)=\int 2^{dj} \rho(2^j(x-y)) d\mu(x) \leq C_{\epsilon} 2^{-j(s-\epsilon)}$$ for any $\epsilon>0$ since $\mu$ is a Frostman measure on $E$. By Schur's test,
$$ U_j: L^2(\mu) \to L^2(\mu) \ \text{with norm} \ C_{\epsilon} 2^{j(d-s-\epsilon)}.$$

By Cauchy-Schwarz,
$$ <U_jf,f> \leq {||U_jf||}_{L^2(\mu)} \cdot {||f||}_{L^2(\mu)} \leq C_{\epsilon} 2^{j(d-s-\epsilon)} {||f||}_{L^2(\mu)}^2$$ and the proof is completed by taking square roots.

\vskip.25in

\vskip.25in 

%\newpage


\begin{thebibliography}{}


\bibitem{BIT15} M. Bennett, A. Iosevich and K. Taylor {\it Finite chains inside thin subsets of ${\mathbb R}^d$}, Anal. PDE \textbf{9} (2016), no. 3, 579-614.



\bibitem{DIOWZ21} X. Du, A. Iosevich, Y. Ou, H. Wang, and R. Zhang, {\it An improved result for Falconer's distance set problem in even dimensions}, Math. Ann. 380 (2021), no. 3-4, 1215-1231.


\bibitem{DORZ23} X. Du, Y. Ou, K. Ren, and R. Zhang, {\it New improvement to Falconer distance set problem in higher dimensions}, (2023), (arXiv:2309.04103). 

\bibitem{Erd05} B. Erdo\~{g}an {\it A bilinear Fourier extension theorem and applications to the distance set problem} IMRN (2006).

\bibitem{EIT11} S. Eswarathasan, A. Iosevich and K. Taylor, {\it Fourier integral operators, fractal sets, and the regular value theorem}, (English summary) Adv. Math. \textbf{228} (2011), no. 4, 2385-2402.

\bibitem{Fal86} K. J. Falconer, {\it On the Hausdorff dimensions of distance sets}, Mathematika \textbf{32} (1986), 206-212.

\bibitem{GGPP23} J. Gaitan, A. Greenleaf, E. Palsson and G. Psaromiligkos,
{\it On restricted Falconer distance sets},
arXiv:2305.18053v2 (July 2023), {\it Canadian Jour. Math.}, 
doi:10.4153/S0008414X24000117.

\bibitem{Graf14} L. Grafakos,  {\it Modern Fourier analysis}, (3rd ed.) Springer, New York, 2014. 

\bibitem{GGIP12} L. Grafakos, A. Greenleaf, A. Iosevich and E. Palsson, {\it Multilinear generalized Radon transforms and point configurations}, Forum Math. {\bf 27} (2015), no. 4, 2323--2360. 

\bibitem{GIT19} A. Greenleaf, A. Iosevich and K. Taylor, {\it Configuration sets with nonempty interior}, 
J. Geometric Analysis {\bf 31} (2021), no. 7, 6662-6680; 
doi:10.1007/s12220-019-00288-y.

\bibitem{GIT20} A. Greenleaf, A. Iosevich and K. Taylor, {\it On $k$-point configuration sets with nonempty 
interior}, arXiv:2005.10796; Mathematika {\bf 68} (2022), no. 1, 163--190; 
doi:10.1112/mtk.12114.

\bibitem{GIT22} A. Greenleaf, A. Iosevich and K. Taylor, {\it Non-empty interior of configuration sets via microlocal partition optimization}, Math. Zeitschrift\, {\bf 306} (2024), no. 4,  doi:10.1007/s00209-024-03466-z.

\bibitem{GS94} A. Greenleaf and A. Seeger, {\it Fourier integral operators with fold singularities}, J. Reine Angew. Math. 455 (1994), 35-56.

\bibitem{GIOW20} L. Guth, A. Iosevich, Y. Ou and H. Wang, {\it On Falconer's distance problem in the plane}, arXiv:1808.09346 Invent. Math. 219 (2020), no. 3, 779-830.

\bibitem{IKSTU19} A. Iosevich, B. Krause, E. Sawyer, K. Taylor and I. Uriarte-Tuero, {\it Maximal operators: scales, curvature and the fractal dimension},  Anal. Math. 45 (2019), no. 1, 63-86.

\bibitem{ILX2020} A. Iosevich, B. Liu and Y. Xi, {\it Microlocal decoupling inequalities and the distance problem on Riemannian manifolds}, (accepted for publication by the Amer. J. Math.), (arXiv:1909.05171), (2020).


\bibitem{IMT11} A. Iosevich, M. Mourgoglou and K. Taylor, 
{\it On the Mattila-Sj\"olin theorem for distance sets}, Ann. Acad. Sci. Fenn. Math. {\bf 37} (2012), no. 2, 557-562.

\bibitem{IT19} A. Iosevich and K. Taylor, {\it Finite trees inside thin subsets of ${\Bbb R}^d$}, 
Modern Methods in Operator Theory and Harmonic Analysis,  Springer Proc. Math. Stat.
{\bf 291} (2019), 51-56.

\bibitem{ITU16} A. Iosevich, K. Taylor and I. Uriarte-Tuero, {\it Pinned geometric configurations in Euclidean space and Riemannian manifolds}, Mathematics, \textbf{9}, 
arXiv:1610.00349 (2021).

\bibitem{KPS21} D. Koh, T. Pham and C.-Y. Shen, {\it On the Mattila-Sj\"{o}lin distance theorem for product sets}, 
arXiv:2103.11418 (March, 2021).

\bibitem{M95} P. Mattila, {\it Geometry of sets and measures in Euclidean spaces}, Cambridge Studies in Adv. Math. {\bf 44}.  Cambridge Univ. Press,1995.

\bibitem{Mat15} P. Mattila, {\it Fourier analysis and Hausdorff dimension.} Cambridge Studies in Adv. 
Math. {\bf 150}.  Cambridge Univ. Press, 2015.


\bibitem{MS99} P. Mattila and P. Sj\"{o}lin, {\it Regularity of distance measures and sets}, Math. Nachr. \textbf{204} (1999), 157-162.

%\bibitem{Mat2015} P. Mattila, {\it Fourier Analysis and Hausdorff Dimension}, Cambridge University Press, (2015).


\bibitem{McTBaire} A. McDonald and K. Taylor, {\it Point configurations in sets of sufficient topological structure and a topological {E}rdős similarity conjecture}, arXiv:2502.10204 (February 2025).

\bibitem{McTconstantGap} A. McDonald and K. Taylor, Infinite constant gap length trees in products of thick Cantor sets, Proc. Roy. Soc. Edinburgh Sect. A {\bf 154} (2024), no.~5, 1336--1347; MR4806278

\bibitem{McTchains} A. McDonald and K. Taylor, Finite point configurations in products of thick Cantor sets and a robust nonlinear Newhouse gap lemma, Math. Proc. Cambridge Philos. Soc. {\bf 175} (2023), no.~2, 285--301; MR4623515


\bibitem{O15} T. Orponen, {\it On the distance sets of AD-regular sets}, Adv. Math. 307, (2017), 1029-1045.

\bibitem{PRA21} E. Palsson and F. Romero Acosta, {\it A Mattila-Sj\"{o}lin theorem for triangles}, 
Jour. Functional Analysis {\bf 284} (6) (2023), 109814.

\bibitem{PRA22} E. Palsson and F. Romero Acosta, {\it  A Mattila-Sj\"olin theorem for simplices in low dimensions},
https://arxiv.org/abs/2208.07198 (2022); Math. Ann., to appear.


\bibitem{T11} M. Taylor, {\it Partial differential equations I. Basic theory}, (2nd ed.), 115. Springer, (2011).

\bibitem{W04} T. Wolff, {\it Lectures on harmonic analysis}, I. Laba and C. Shubin, eds. University Lec. Series, \textbf{29}. American Mathematical Society, Providence, RI, 2003.

\end{thebibliography}
\end{document}